\documentstyle[11pt,twoside]{article}
\include{graphicx}
\title{An asymptotic expansion for the Stieltjes constants}
\author{\sc R. B.\ Paris \\
{\em Division of Computing and Mathematics}, \\
{\em University of Abertay Dundee, Dundee DD1 1HG, UK}
}
\begin{document}
\def\f#1#2{\mbox{${\textstyle \frac{#1}{#2}}$}}
\def\dfrac#1#2{\displaystyle{\frac{#1}{#2}}}
\def\boldal{\mbox{\boldmath $\alpha$}}
{\newcommand{\Sgoth}{S\;\!\!\!\!\!/}
\newcommand{\bee}{\begin{equation}}
\newcommand{\ee}{\end{equation}}
\newcommand{\lam}{\lambda}
\newcommand{\ka}{\kappa}
\newcommand{\al}{\alpha}
\newcommand{\fr}{\frac{1}{2}}
\newcommand{\fs}{\f{1}{2}}
\newcommand{\g}{\Gamma}
\newcommand{\br}{\biggr}
\newcommand{\bl}{\biggl}
\newcommand{\ra}{\rightarrow}
\newcommand{\mbint}{\frac{1}{2\pi i}\int_{c-\infty i}^{c+\infty i}}
\newcommand{\mbcint}{\frac{1}{2\pi i}\int_C}
\newcommand{\mboint}{\frac{1}{2\pi i}\int_{-\infty i}^{\infty i}}
\newcommand{\gtwid}{\raisebox{-.8ex}{\mbox{$\stackrel{\textstyle >}{\sim}$}}}
\newcommand{\ltwid}{\raisebox{-.8ex}{\mbox{$\stackrel{\textstyle <}{\sim}$}}}
\renewcommand{\topfraction}{0.9}
\renewcommand{\bottomfraction}{0.9}
\renewcommand{\textfraction}{0.05}
\newcommand{\mcol}{\multicolumn}
\date{}
\maketitle
\pagestyle{myheadings}
\markboth{\hfill \sc R. B.\ Paris  \hfill}
{\hfill \sc  Asymptotics of the Stieltjes constants\hfill}
\begin{abstract}
The Stieltjes constants $\gamma_n$ appear in the coefficients in the Laurent expansion of the Riemann zeta function $\zeta(s)$ about the simple pole $s=1$. We present an asymptotic expansion for $\gamma_n$ as $n\rightarrow \infty$ based on the approach described by Knessl and Coffey [Math. Comput. {\bf 80} (2011) 379--386]. A truncated form of this expansion with explicit coefficients is also given. Numerical results are presented that illustrate the accuracy
achievable with our expansion.
\vspace{0.4cm}

\vspace{0.3cm}

\end{abstract}

\vspace{0.3cm}

\noindent $\,$\hrulefill $\,$

\vspace{0.2cm}

\begin{center}
{\bf 1. \  Introduction}
\end{center}
\setcounter{section}{1}
\setcounter{equation}{0}
\renewcommand{\theequation}{\arabic{section}.\arabic{equation}}
The Stieltjes constants $\gamma_n$ appear in the coefficients in the Laurent expansion of the Riemann zeta function $\zeta(s)$ about the point $s=1$ given by
\[\zeta(s)=\frac{1}{s-1}+\sum_{n=0}^\infty \frac{(-)^n}{n!}\,\gamma_n\,(s-1)^n,\]
where $\gamma_0=0.577216 \ldots$ is the well-known Euler-Mascheroni constant. Some historical notes and numerical values of $\gamma_n$ for $n\leq 20$ are given in \cite{BF}. Recent high-precision evaluations of $\gamma_n$ based on numerical integration have been described in \cite{Ke, K}. In \cite{Ke}, Keiper lists various $\gamma_n$ up to $n=150$, whereas in \cite{K}, Kreminski has computed values to several thousand digits for $n\leq 10^4$ and for further selected values (accurate to $10^3$ digits) up to $n=5\times 10^4$. All values up to $n=10^5$ have been computed by Johansson in \cite{J} to about $10^4$ digits.

Upper bounds for $|\gamma_n|$ in the form
\[|\gamma_n|\leq \{3+(-)^n\}\,\frac{\lambda_n\,\g(n)}{\pi^n},\]
have been obtained by Berndt \cite{B}
with $\lambda_n=1$, and by Zhang and Williams \cite{ZW} with $\lambda_n=(2/n)^n \pi^{-\fr}\g(n+\fs)\sim \surd 2 (2/e)^n$ for $n\ra\infty$. On the other hand, Matsuoka \cite{M} has shown that
\[|\gamma_n|\leq 10^{-4} e^{n\log\,\log\,n}\qquad (n\geq 10).\]
However, all these bounds grossly overestimate the growth of $|\gamma_n|$ for large values of $n$.
An asymptotic approximation for $\gamma_n$ has recently been given by Knessl and Coffey \cite{KC} in the form
\bee\label{e11}
\gamma_n\sim \frac{Be^{nA}}{\sqrt{n}} \,\cos\,(na+b)\qquad(n\ra+\infty),
\ee
where $A$, $B$, $a$ and $b$ are functions that depend weakly on $n$; see Section 2 for the definition of these quantities. Knessl and Coffey  have verified numerically that for $n\leq 3.5\times 10^4$ the above formula accounts for the asymptotic growth and oscillatory pattern of $\gamma_n$, with the exception of $n=137$ where the cosine factor in (\ref{e11}) becomes very small.

The aim in this note is to extend the analysis in \cite{KC} to generate an asymptotic expansion for $\gamma_n$ as $n\ra+\infty$. The coefficients in this expansion are determined numerically by application of Wojdylo's formulation \cite{W} for the coefficients in the expansion of a Laplace-type integral. An explicit evaluation of the coefficients is obtained in the case of the expansion truncated after three terms. This approximation is extended to the more general Stieltjes constants $\gamma_n(\alpha)$ appearing in the Laurent expansion of the Hurwitz zeta function $\zeta(s,\alpha)$. Numerical results are presented in Section 3 to demonstrate the accuracy of our expansion compared to that in (\ref{e11}).

\vspace{0.6cm}

\begin{center}
{\bf 2. \ Asymptotic expansion for $\gamma_n$}
\end{center}
\setcounter{section}{2}
\setcounter{equation}{0}
\renewcommand{\theequation}{\arabic{section}.\arabic{equation}}
We start with the integral representation for $n\geq 1$ given in \cite{ZW}
\[\gamma_n=\int_1^\infty B_1(x-[x])\,\frac{\log^{n-1}x}{x^2}\,(n-\log\,x)\,dx,\]
where $B_1(x-[x])=-\sum_{j=1}^\infty \frac{\sin 2\pi jx}{\pi j}$ is the first periodic Bernoulli polynomial. With
the change of variable $t=\log\,x$, we obtain \cite[Eq.~(2.3)]{KC}
\[\gamma_n=-\Im\,\bl\{\sum_{k=1}^\infty \frac{1}{\pi k} \int_0^\infty \exp\,[2\pi ike^t+n\log\,t-t]\,\bl(\frac{n}{t}-1\br) dt\br\}.\]
Following the approach used in \cite{KC}, we define 
\bee\label{e21a}
\psi_k(t)\equiv\psi_k(t;n)=-\frac{2\pi ik}{n}\,e^t-\log\,t,\qquad f(t)\equiv f(t;n)=\frac{e^{-t}}{t}\bl(1-\frac{t}{n}\br)
\ee
and write
\bee\label{e21}
\gamma_n=-\Im \sum_{k=1}^\infty J_k,\qquad J_k:=\frac{n}{\pi k}\int_0^\infty e^{-n\psi_k(t)}f(t)\,dt.
\ee

We employ the method of steepest descents to estimate the integrals $J_k$ for large values of $n$. 
Saddle points of the exponential factor occur at the zeros of $\psi_k'(t)=0$; that is, they satisfy
\bee\label{e22}
te^t=\frac{ni}{2\pi k}.
\ee
There is an infinite string of saddle points, which is approximately parallel to the imaginary $t$-axis, given by \cite{KC}
\[t_m=\log\,\frac{n}{2\pi k}-\log\log\,n+(2m+\fs)\pi i+O\bl(\frac{\log\log\,n}{\log\,n}\br)\]
for $m=0, \pm 1, \pm 2, \ldots$ and large $n$. For fixed $k$ and $m$, the value of $\Re\,\psi_k(t_m)$ is then 
\[-\Re\,\psi_k(t_m)=\log\log\,n-\frac{1}{\log\,n}(1+\log\,(2\pi k\log\,n))+O((\log\,n)^{-2})\]
as $n\ra\infty$, where the dependence on $m$ is contained in the order term.
This shows that the heights of the saddles corresponding to $k\geq 2$ are exponentially smaller as $n\ra\infty$ than the saddle
with $k=1$, so that to within exponentially small correction terms we may neglect the contribution in (\ref{e21})
arising from $k$ values corresponding to $k\geq 2$; but see the discussion in Section 3. From hereon, we shall drop the subscript $k$ and write $\psi_1(t)\equiv\psi(t)$. 

Typical paths of steepest descent and ascent through the saddles $t_0$ and $t_1$ are shown in Fig.~1. Steepest descent and ascent paths terminate at infinity in the right-half plane in the directions $\Im (t)=(2j+\fs)\pi$ and $\Im (t)=(2j+\f{3}{2})\pi$ ($j=0, \pm 1, \pm 2,\ldots$), respectively.
The steepest descent paths through $t_0$ and $t_1$ emanate from the origin and pass to infinity in the directions $\Im (t)=\fs\pi$ and $\f{5}{2}\pi$, respectively. Similarly, the steepest descent path through $t_{-1}$ (not shown) emanates from the origin and passes to infinity in the direction $\Im (t)=-\f{3}{2}\pi$. The integration path in (\ref{e21}) can then be deformed to pass through the saddle $t_0$ as shown in Fig.~1.
\begin{figure}[t]
	\begin{center}\includegraphics[width=0.4\textwidth]{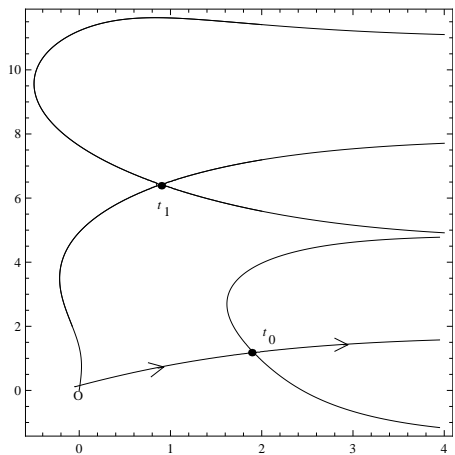}\\
\caption{\small{Paths of steepest descent and ascent through the saddles $t_0$ and $t_1$ when $n=100$ and $k=1$. 
The steepest paths through the saddle $t_{-1}$ (not shown) in the lower half-plane are similar to those through $t_1$. The arrows indicate the direction of integration.}}
	\end{center}
\end{figure}

Application of the method of steepest descents (see, for example, \cite[p.~127]{O} and \cite[p.~14]{P}) then yields
\bee\label{e23}
J_1\sim \frac{n}{\pi}\,\frac{\sqrt{2\pi}\, e^{-n\psi(t_0)-t_0}}{t_0(-\psi''(t_0))^{1/2}}\,\bl(1-\frac{t_0}{n}\br) \sum_{s=0}^\infty \frac{{\hat c}_{2s} (\fs)_s}{n^{s+1/2}},
\ee
where $(a)_s=\g(a+s)/\g(a)$ is the Pochhammer symbol and ${\hat c}_0=1$.
The normalised coefficients ${\hat c}_{2s}$ can be obtained by an inversion process and are listed for $s\leq 4$ in \cite[p.~119]{D} and for $s\leq 2$ in \cite[p.~13]{P}; see below.
Alternatively, they can be obtained by an expansion process to yield Wojdylo's formula \cite{W} given by
\bee\label{e25}
{\hat c}_{s}=\alpha_0^{-s/2}\sum_{k=0}^{s} \frac{\beta_{s-k}}{\beta_0} \sum_{j=0}^k \frac{(-)^j (\fs s+\fs)_j}{j!\ \alpha_0^j}\,{\cal B}_{kj}\,;
\ee
see also \cite[p.~25]{T}. Here ${\cal B}_{kj}\equiv {\cal B}_{kj}(\alpha_1, \alpha_2, \ldots , \alpha_{k-j+1})$ are the partial ordinary Bell polynomials generated by the recursion\footnote{For example, this generates the values ${\cal B}_{41}=\alpha_4$, ${\cal B}_{42}=\alpha_3^2+2\alpha_1\alpha_3$, ${\cal B}_{43}=3\alpha_1^2\alpha_2$ and ${\cal B}_{44}=\alpha_1^4$.}
\[{\cal B}_{kj}=\sum_{r=1}^{k-j+1} \alpha_r {\cal B}_{k-r,j-1} ,\qquad {\cal B}_{k0}=\delta_{k0},\]
where $\delta_{mn}$ is the Kronecker symbol, and the coefficients $\alpha_r$ and $\beta_r$ appear in the expansions
\bee\label{e260}
\psi(t)-\psi(t_0)=\sum_{r=0}^\infty \alpha_r (t-t_0)^{r+2},\qquad f(t)=\sum_{r=0}^\infty \beta_r(t-t_0)^r
\ee
valid in a neighbourhood of the saddle $t=t_0$.

Following \cite{KC}, we put $t_0=u+iv$, where $u$, $v$ are real, and write $-\psi(t_0)=A+ia$, where
\bee\label{e250}\left.
\begin{array}{ll}A:&=\Re (\log\,t_0-1/t_0)=\frac{1}{2}\log (u^2+v^2)-\dfrac{u}{u^2+v^2},\\
a:&=\Im (\log\,t_0-1/t_0)=\arctan \bl(\dfrac{v}{u}\br)+\dfrac{v}{u^2+v^2}.\end{array}\right\}
\ee
We have $\psi''(t_0)=(1+t_0)/t_0^2$ and accordingly define\footnote{In \cite{KC}, the quantity $\fs\pi-v$ appearing in the definition of $b$ is written as $\arctan\,(v/u)$ by virtue of the first relation in (\ref{e220}).} 
\bee\label{e26}
B:=2\sqrt{2\pi}\,\bl|\frac{t_0}{\sqrt{1+t_0}}\br|,\qquad b:=\fs\pi-v-\arctan\,\bl(\frac{v}{1+u}\br).
\ee 
A simple calculation using (\ref{e22}) with $k=1$ shows that
\bee\label{e220}
\tan v=\frac{u}{v},\qquad e^{-u}=\frac{2\pi |t_0|}{n}.
\ee
Then, from (\ref{e21}) with $k=1$, (\ref{e23}) and the second relation in (\ref{e220}),
we find upon incorporating the factor $1-t_0/n$ into the asymptotic series that
\[\gamma_n\sim \frac{Be^{nA}}{\sqrt{n}}\,\Re \bl\{e^{i(na+b)} \sum_{s=0}^\infty \frac{c'_{2s} (\fs)_s}{n^s}\br\},\]
where
\bee\label{e27a}
c'_{2s}={\hat c}_{2s}-\frac{2t_0}{2s-1}\,{\hat c}_{2s-2}\qquad (s\geq 1).
\ee

If we now introduce the real and imaginary parts of the coefficients ${\hat c}_{2s}$ by
\bee\label{e29}
c'_{2s}:=C_s+iD_s\qquad (s\geq 1),\qquad C_0=1,\ \ D_0=0,
\ee
where we recall that $C_s$ and $D_s$ contain an $n$-dependence,
then we have the expansion of $\gamma_n$ given by the following theorem. 
\newtheorem{theorem}{Theorem}
\begin{theorem}$\!\!\!.$\ \ Let the quantities $A$, $B$, $a$ and $b$, and the coefficients $C_s$, $D_s$, be as defined in (\ref{e250}), (\ref{e26}) and (\ref{e29}).
Then, neglecting exponentially smaller terms, we have
\bee\label{e27}
\gamma_n\sim\frac{Be^{nA}}{\sqrt{n}}\,\bl\{\cos\,(na+b) \sum_{s=0}^\infty \frac{C_s (\fs)_s}{n^s}-\sin\,(na+b) \sum_{s=1}^\infty \frac{D_s (\fs)_s}{n^s}\br\}
\ee
as $n\ra\infty$.
\end{theorem}
We note that to leading order $A\sim \log\log\,n$ and $B\sim(8\pi\log\,n)^{1/2}$ for large $n$.

A simpler form of the expansion (\ref{e27}) can be given by truncating the above series at $s=2$ and use of the form of the normalised coefficients ${\hat c}_{2s}$ in (\ref{e23}) expressed in the form 
\[{\hat c}_2=\frac{1}{2\psi''(t_0)}\{2F_2-2\Psi_3F_1+\f{5}{6}\Psi_3^2-\fs\Psi_4\},\]
\[{\hat c}_4=\frac{1}{(2\psi''(t_0))^2}\{\f{2}{3}F_4-\f{20}{9}\Psi_3F_3+\f{5}{3}(\f{7}{3}\Psi_3^2-\Psi_4)F_2-\f{35}{9}(\Psi_3^3-\Psi_3\Psi_4+\f{6}{35}\Psi_5)F_1 \]
\[+\f{35}{9}(\f{11}{24}\Psi_3^4-\f{3}{4}(\Psi_3^2-\f{1}{6}\Psi_4)\Psi_4+\f{1}{5}\Psi_3\Psi_5-\f{1}{35}\Psi_6)\}\]
where, for brevity, we have defined
\[\Psi_m:=\frac{\psi^{(m)}(t_0)}{\psi''(t_0)}\ \ \ (m\geq 3),\qquad F_m:=\frac{f^{(m)}(t_0)}{f(t_0)}\ \ \ (m\geq 1);\]
see \cite[p.~119]{D}, \cite[pp.~13--14]{P}.

From (\ref{e21a}) and (\ref{e27a}), use of {\it Mathematica} shows that
\[c'_2=\frac{\wp_2(t_0)}{12(1+t_0)^3}+\frac{(4+3t_0)t_0^2}{n(1+t_0)^2}+O(n^{-2}),\qquad  c'_4=\frac{\wp_4(t_0)}{864(1+t_0)^6}+O(n^{-1}),\]
where
\[\wp_2(t_0)=2-18t_0-20t_0^2-3t_0^3+2t_0^4,\]
\[\wp_4(t_0)=4-72t_0-332t_0^2-8028t_0^3-19644t_0^4-20280t_0^5-9911t_0^6-1884t_0^7+4t_0^8.\]
Then we obtain the following result.
\begin{theorem}$\!\!\!.$\ \ Let the quantities $A$, $B$, $a$ and $b$ be as defined in (\ref{e250}) and (\ref{e26}).
Then, with
\[{c}_1+i{d}_1=\frac{\wp_2(t_0)}{24(1+t_0)^3},\qquad {c}_2+i{d}_2=\frac{\wp_4(t_0)}{1152(1+t_0)^6}+\frac{(4+3t_0)t_0^2}{2(1+t_0)^2},\]
where ${c}_s$, ${d}_s$ $(s=1, 2)$ are real (and independent of $n$) and $t_0$ is the saddle point given by the principal solution of (\ref{e22}) with $k=1$, we have the asymptotic approximation
\bee\label{e210}
\gamma_n\sim \frac{Be^{nA}}{\sqrt{n}}\,\bl\{\cos\,(na+b)\bl(1+\frac{{c}_1}{n}+\frac{{c}_2}{n^2}\br)
-\sin\,(na+b)\bl(\frac{{d}_1}{n}+\frac{{d}_2}{n^2}\br)\br\}
\ee
as $n\ra\infty$.
\end{theorem}

We remark that the expansion of the integrals $J_k$ for fixed $k\geq 2$ follows the same procedure. If we still refer to the real and imaginary parts of the contributory saddle $t_0$ (when $k\geq 2$) as $u$ and $v$, the second relation in (\ref{e220}) is now replaced by $e^{-u}=2\pi k |t_0|/n$. It then follows that the form of the expansion for $-\Im J_k$ is given by (\ref{e27}), provided the quantities $A$, $B$, $a$ and $b$, and the coefficients $C_s$, $D_s$, are interpreted as corresponding to the saddle $t_0$ with the $k$-value under consideration.

\vspace{0.6cm}

\begin{center}
{\bf 3. \ Numerical results and concluding remarks}
\end{center}
\setcounter{section}{3}
\setcounter{equation}{0}
\renewcommand{\theequation}{\arabic{section}.\arabic{equation}}
We discuss numerical computations carried out using the expansions given in Theorems 1 and 2.
For a given value of $n$ the saddle $t_0$ is computed from (\ref{e22}) with $k=1$ to the desired accuracy. {\it Mathematica} is used to determine the coefficients $\alpha_r$ and $\beta_r$ in (\ref{e260}) for $0\leq r\leq 2s_0$, where in the present computations $s_0=6$. The coefficients $C_s$ and $D_s$ can then be calculated for $0\leq s\leq s_0$ from (\ref{e25}), (\ref{e27a}) and (\ref{e29}). 

We display the computed values of $C_s$ and $D_s$ for two values of $n$ in Table 1. 
We repeat that these coefficients contain an $n$-dependence and so have to be computed for each value of $n$ chosen.
In Table 2, the values of the absolute relative error in the computation of $\gamma_n$ from the expansion (\ref{e27}) are presented as a function of the truncation index $s$ for several values of $n$. 
\begin{table}[th]
\caption{\footnotesize{Values of the coefficients $C_s$ and $D_s$ (to 10 dp) for $1\leq s\leq 6$ for two values of $n$.}}
\begin{center}
\begin{tabular}{|l|rr|rr|}
\hline
&&&&\\[-0.3cm]
\mcol{1}{|c|}{} & \mcol{2}{|c|}{$n=100$} & \mcol{2}{c|}{$n=1000$}\\  
\mcol{1}{|c|}{$s$}& \mcol{1}{c}{$C_s$} & \mcol{1}{c|}{$D_s$} & \mcol{1}{c}{$C_s$} & \mcol{1}{c|}{$D_s$} \\
[.1cm]\hline
&&&&\\[-0.3cm]
1 & $-0.3158578918$ & $+0.1626819326$ & $-0.0885061806$ & $+0.1958085240$\\
2 & $-2.9096870797$ & $-2.1947177121$ & $-6.5840165991$ & $-2.6459812815$\\
3 & $-0.3804847598$ & $-3.3953890569$ & $-9.4682639154$ & $-10.09635962642$\\
4 & $+1.4820479884$ & $-0.1130053628$ & $-1.3074432243$ & $-11.31040992292$\\
5 & $-0.2630549338$ & $+0.9253656779$ & $+4.9469591967$ & $-1.67819725309$\\
6 & $-0.3783700609$ & $-0.3119889058$ & $+0.8180579543$ & $+3.98701271605$\\
[.2cm]\hline
\end{tabular}
\end{center}
\end{table}
\begin{table}[th]
\caption{\footnotesize{Values of the absolute relative error in the computation of $\gamma_n$ from (\ref{e27}) 
as a function of the truncation index $s$ for different $n$.}}
\begin{center}
\begin{tabular}{|l|l|l|l|l|}
\hline
&&&&\\[-0.3cm]
%\mcol{1}{|c|}{} & \mcol{2}{|c|}{$n=100$} & \mcol{2}{c|}{$n=1000$}\\  
\mcol{1}{|c|}{$s$} & \mcol{1}{c|}{$n=75$}& \mcol{1}{c|}{$n=100$} & \mcol{1}{c|}{$n=137$} & \mcol{1}{c|}{$n=1000$} \\
[.1cm]\hline
&&&&\\[-0.3cm]
0 & $1.759\times 10^{-3}$ & $1.412\times 10^{-3}$  & $\ \ \ \ \ \ -\!\!-$   & $1.597\times 10^{-4}$ \\
1 & $6.503\times 10^{-4}$ & $3.226\times 10^{-4}$  & $2.701\times 10^{-1}$  & $2.649\times 10^{-6}$\\
2 & $1.244\times 10^{-5}$ & $4.472\times 10^{-6}$  & $8.775\times 10^{-2}$  & $4.125\times 10^{-9}$\\
3 & $3.063\times 10^{-7}$ & $9.370\times 10^{-8}$  & $3.811\times 10^{-5}$  & $7.711\times 10^{-11}$\\
4 & $2.535\times 10^{-9}$ & $7.850\times 10^{-10}$ & $2.183\times 10^{-6}$  & $2.026\times 10^{-13}$\\
5 & $5.101\times 10^{-10}$& $9.022\times 10^{-11}$ & $1.248\times 10^{-8}$  & $6.157\times 10^{-16}$\\
6 & $1.850\times 10^{-11}$& $1.982\times 10^{-12}$ & $9.415\times 10^{-10}$ & $2.743\times 10^{-18}$\\
[.2cm]\hline
\end{tabular}
\end{center}
\end{table}

The case $n=137$ has been included in Table 2 since this corresponds to the factor $\cos (na+b)$ possessing the very small value $\simeq 1.69881\times 10^{-4}$. The leading term approximation in (\ref{e11}), and (\ref{e27}) (with $s=0$), yields an incorrect sign, namely $+3.89874\times 10^{27}$ when $\gamma_{137}=-7.99522199\ldots \times 10^{27}$. According to \cite{J}, this is the only instance for $n\leq 10^5$ when the leading approximation produces the wrong sign. It is seen that inclusion of the higher order correction terms with $s\leq 6$ yields a relative error of order $10^{-10}$ in this case.
When $n=10^5$, \cite{J} gives the value
\[\gamma_{100000}=1.99192730631254109565822724315 \ldots \times 10^{83432}.\]
The expansion (\ref{e27}) for this value of $n$ with truncation index $s=6$ is found to yield a relative error of order $10^{-30}$; that is, the expansion correctly reproduces all the digits displayed above.

\begin{table}[th]
\caption{\footnotesize{Values of the absolute relative error in the computation of $\gamma_n$ from (\ref{e27}) with $k=1$ and $k\leq 2$ as a function of the truncation index $s$ for $n=25$.}}
\begin{center}
\begin{tabular}{|l|l|l|}
\hline
&&\\[-0.3cm]
%\mcol{1}{|c|}{} & \mcol{2}{|c|}{$n=100$} & \mcol{2}{c|}{$n=1000$}\\  
\mcol{1}{|c|}{$s$} & \mcol{1}{c|}{$k=1$}& \mcol{1}{c|}{$k\leq 2$} \\
[.1cm]\hline
&&\\[-0.3cm]
0 & $1.051\times 10^{-2}$ & $1.052\times 10^{-2}$ \\
1 & $2.909\times 10^{-3}$ & $2.894\times 10^{-3}$ \\
2 & $2.608\times 10^{-4}$ & $2.460\times 10^{-4}$ \\
3 & $2.390\times 10^{-6}$ & $1.723\times 10^{-5}$ \\
4 & $1.518\times 10^{-5}$ & $3.412\times 10^{-7}$ \\
5 & $1.495\times 10^{-5}$ & $1.160\times 10^{-7}$ \\
6 & $1.482\times 10^{-5}$ & $1.189\times 10^{-8}$ \\
[.2cm]\hline
\end{tabular}
\end{center}
\end{table}

For the smallest value $n=75$ presented in Table 2, it is found numerically that the contribution to 
(\ref{e21}) corresponding to $k=2$ is about 11 orders of magnitude smaller than the dominant term with $k=1$. For the larger $n$ values, this contribution is even smaller and the terms with $k\geq 2$ can be safely neglected.
However, for smaller $n$ this is no longer the case and a meaningful approximation has to take into account the contribution from other $k\geq 2$ values. 

In Table 3, we illustrate this situation by taking $n=25$. The second column shows the absolute relative error in the computation of $\gamma_n$ with $k=1$ for different truncation index $s$; that is, with the approximation $\gamma_n\simeq -\Im J_1$.
For $4\leq s\leq 6$, this error is seen to remain essentially constant at $O(10^{-5})$. The contribution with $k=2$ is about 5 orders of magnitude smaller than the $k=1$ contribution, so that this additional contribution needs to be included for larger index $s$.
The absolute relative error including the contribution with $k=2$ is shown in the third column; that is, with the approximation $\gamma_n\simeq -\Im (J_1+J_2)$. 
The expansion with $k=3$ is about 8 orders of magnitude smaller than the $k=1$ contribution, so this would
only begin to make a significant contribution for $s\geq 6$. This problem becomes even more acute for smaller $n$
values, say $n=10$, where higher $k$ values need to be retained. However, the chief interest in the asymptotic expansion in (\ref{e27}) is for large $n$, where this problem is of no real concern.

In Table 4 we show some examples of the asymptotic approximation given in (\ref{e210}). We compare these with the values produced by the leading approximation (\ref{e11}) and the exact value of $\gamma_n$ obtained from {\it Mathematica} using the command {\tt StieltjesGamma[n]}. It will be observed that for $n=500$ the approximation (\ref{e210}) yields nine significant figures.
\begin{table}[t]
\caption{\footnotesize{Values for $\gamma_n$ obtained from (\ref{e11}) and (\ref{e210}) 
compared with the exact value.}}
\begin{center}
\begin{tabular}{|l|l|l|l|}
\hline
&&&\\[-0.3cm]
%\mcol{1}{|c|}{} & \mcol{2}{|c|}{$n=100$} & \mcol{2}{c|}{$n=1000$}\\  
\mcol{1}{|c|}{$n$} & \mcol{1}{c|}{Eq.~(1.1)} & \mcol{1}{c|}{Eq.~(2.13)} & \mcol{1}{c|}{Exact\ $\gamma_n$} \\
[.1cm]\hline
&&&\\[-0.3cm]
10 & $+2.105395\times 10^{-4}$ & $+2.04713213\times 10^{-4}$ & $+2.05332814\ldots\times 10^{-4}$\\
50 & $+1.275493\times 10^{2}$ & $+1.26823798\times 10^{2}$ & $+1.26823602\ldots\times 10^{2}$\\
80 & $+2.514857\times10^{10}$ & $+2.51633995\times10^{10}$ & $+2.51634410\ldots\times10^{10}$\\
100& $-4.259408\times10^{17}$ & $-4.25340036\times10^{17}$ & $-4.25340157\ldots\times10^{17}$\\
137& $+3.898740\times10^{27}$ & $-7.99377883\times10^{27}$ & $-7.99522199\ldots\times10^{27}$\\
200& $-7.060244\times 10^{55}$ & $-6.97465335\times 10^{55}$ & $-6.97464971\ldots\times 10^{55}$\\ 
500& $-1.165662\times 10^{204}$ & $-1.16550527\times10^{204}$ & $-1.16550527\ldots\times 10^{204}$\\
[.2cm]\hline
\end{tabular}
\end{center}
\end{table}

Finally, we remark that the analysis in Section 2 is immediately applicable to the more general Stieltjes constants $\gamma_n(\alpha)$ appearing in the Laurent expansion for the Hurwitz zeta function $\zeta(s,\alpha)$ about the point $s=1$. These constants are defined by 
\[\zeta(s,\alpha)=\frac{1}{s-1}+\sum_{n=0}^\infty \frac{(-)^n}{n!}\,\gamma_n(\alpha)\,(s-1)^n,\]
where $\gamma_0(\alpha)=-\g'(\alpha)/\g(\alpha)$ and $\gamma_n(1)=\gamma_n$.
It is shown in \cite[Eq.~(2.9)]{KC1} that
\[C_n(\alpha):=\gamma_n(\alpha)-\frac{1}{\alpha}\,e^{n\log\log\,\alpha}=-\Im \sum_{k=1}^\infty e^{-2\pi ik\alpha} J_k.\]
Then it follows that the expansions in Theorems 1 and 2 are modified only in the argument of the trigonometric functions appearing therein, which become $na+b-2\pi\alpha$. Thus, for example, from (\ref{e210}) we have
\[C_n(\alpha)\sim \frac{Be^{nA}}{\sqrt{n}}\,\bl\{\cos\,(na\!+\!b\!-\!2\pi\alpha)\bl(1+\frac{{c}_1}{n}+\frac{{c}_2}{n^2}\br)
-\sin\,(na\!+\!b\!-\!2\pi\alpha)\bl(\frac{{d}_1}{n}+\frac{{d}_2}{n^2}\br)\br\}\]
as $n\ra\infty$, where the quantities $A$, $B$, $a$, $b$ and the coefficients $c_s, d_s$ ($s=1, 2$) are as specified in Theorem 2.
The leading approximation agrees with that obtained in \cite[Eq.~(2.4)]{KC1}.

\end{document}